\documentstyle[12pt]{article}
  \newcommand{\const}{\rm const}
  \newcommand{\Var}{\rm Var}
  
  \newcommand{\Law}{\rm Law}

  \newcommand{\Dom}{\rm  Dom}
  
  \newcommand{\card}{\rm card}
  \newcommand{\argmin}{\rm argmin}
  \newcommand{\Poisson}{\rm Poisson}

\begin{document}

   \begin{center}

 \ {\bf  Moment and tail estimation for U-statistics  }\\

 \vspace{4mm}

 \  {\bf with positive kernels} \\

 \vspace{7mm}

  {\bf   Ostrovsky E., Sirota L.}\\

\vspace{4mm}

 Israel,  Bar-Ilan University, department of Mathematic and Statistics, 59200, \\

\vspace{4mm}

e-mails: \  eugostrovsky@list.ru \\
sirota3@bezeqint.net \\

\vspace{4mm}

  {\bf Abstract} \\

\vspace{4mm}

 \end{center}

 \ \  We  deduce the non-asymptotical (bilateral) estimates for moment inequalities for multiple
  sums of non-negative (more precisely, non-negative) independent random variables, on the other words, the well known U or V-statistics. Our consideration
based on the correspondent estimates for the one-dimensional case by means of the so-called degenerate approximation. \par
 \ We apply also the theory of Bell functions as well as the properties of the Poisson distribution and the theory of the so-called Grand Lebesgue Spaces (GLS).\par

 \vspace{4mm}

 \par

\vspace{4mm}


\vspace{5mm}

\section{ Definitions.  Notations. Previous results.  Statement of problem.}

\vspace{4mm}

 \ Let $ \   (X, \cal{B},  {\bf \mu})   \ $ and $ \   (Y, \cal{C}, {\bf \nu} ) \  $
be  two non-trivial probability spaces:  $ \ \mu(X)  = \nu(Y) = 1 . \ $ We will denote by $ \  |g|_p  = |g|L(p) \ $ the ordinary
Lebesgue-Riesz $ \  L(p) \  $ norm of arbitrary measurable numerical valued function $ \   g: X \to R: \  $

$$
|g|_p = |g| L(p) = |g| L_p(X, \mu) := \left[ \int_X |g(x)|^p \ \mu(dx)  \right]^{1/p}, \  p \in [1, \infty)
$$
 analogously for the (also measurable) function $ \  h: Y \to R \ $

$$
|h|_p = |h|L(p) = |h|L_p(Y,\nu) := \left[ \int_Y |h(y)|^p \ \nu(dy) \right]^{1/p};
$$
and for arbitrary integrable function {\it of two variables} $ \  f: X \otimes Y  \to R \ $

$$
|f|_p  = |f|L(p) = |f|L_p(X,Y):= \left[ \int_X \int_Y |f(x,y)|^p \ \mu(dx) \ \nu(dy) \right]^{1/p}, \ p \in [1, \infty).
$$

\vspace{4mm}

 \ Let $ \  Z_+ = \{  1,2, 3, \ldots \} \  $ and denote $ \ Z^2_+ = Z_+ \otimes Z_+, \  Z^d_+ = \otimes_{k=1}^d Z_+. \ $
Let also  $ \  \{\xi(i) \} \  $ and $ \  \{\eta(j)\},
 \ i,j = 1,2,\ldots, \   \xi := \xi(1), \ \eta := \eta(1)  $ be {\it common} independent random variables  defined on certain probability space
$ \   (\Omega, \cal{M}, {\bf P})  \ $  with distributions  correspondingly $ \ \mu, \ \nu: \ $

$$
{\bf P}(\xi(i)  \in A) = \mu(A), \ A \in \cal{B};
$$

$$
{\bf P}(\eta(j) \in F) = \nu(F), \ F \in \cal{C}, \eqno(1.0)
$$

so that

$$
{\bf E} |g(\xi)|^p  = |g|_p^p = \int_X |g(x)|^p \mu(dx), \ {\bf E} |h(\eta)|_p^p = |h|^p_p = \int_Y |h(y)|^p \ \nu(dy),
$$
and

$$
{\bf E} |f(\xi, \eta)|^p = |f|^p_p = \int \int_{X\otimes Y} |f(x,y)|^p \  \mu(dx) \ \nu(dy).
$$

 \ Let also $ \  L \ $ be arbitrary non-empty {\it finite} subset of the set $ \  Z^2_+; \ $ denote by $  \   |L| \  $  a  numbers of its
elements (cardinal number): $  \  |L| := \card(L). \  $  It is reasonable  to suppose in what follows $ \ |L| \ge 1. \ $ \par

 \ Define for any {\it  integrable} function $ \  f: X \otimes Y \to R, \  $ i.e. for which

$$
|f|_1 =  {\bf E}|f(\xi, \eta)| = \int_X \int_Y |f(x,y)| \ \mu(dx) \ \nu(dy) < \infty,
$$
the following normalized sum

$$
S_L[f] \stackrel{def}{=} |L|^{-1} \sum_{(k(1), k(2)) \in L} f(\xi(k(1), \eta(k(2)), \eqno(1.1)
$$
which is a slight generalization of the classical $ \ U \ $ and $ \  V \ $ statistics, see the classical  monograph of
Korolyuk V.S and Borovskikh Yu.V. [18]. Offered here report is the direct generalization of a recent article [30], but we
apply here other methods.\par

 \ The {\it  reasonableness } of this norming function $ \   |L|^{-1} \ $
implies that in general, i.e. in the non-degenerate case $ \   \Var(S_L) \asymp 1, \ |L| \ge 1.  \ $ This propositions
holds true still in the multidimensional case. \par

 \ Our notations  and some previous results are borrowed from the works of S.Klesov  [14] - [17], see also [18]. \par

\vspace{4mm}

\ {\it  We will suppose in what follows that the function   }  $ \ f = f(x,y), \ x \in X, \ y \in Y \ $, {\it as well as both the functions} $ \ g, \ h  \ $ {\it are non - negative (a.e). }

\vspace{4mm}

 \ The so-called {\it centered case}, i.e. when

$$
 {\bf E} f(\xi, \eta) = \int_X \int_Y f(x,y) \ \mu(dx) \ \nu(dy) = 0,
$$
was investigated in many works, see e.g. the classical monograph of Korolyuk V.S and Borovskikh Yu.V. [18]; see also [14]-[17], [26]-[27], [28]
etc. The one-dimensional case was studies in the recent preprint [30]. \par

\vspace{4mm}

 \ {\bf   Our claim in this report is to derive  the moment and exponential bounds for tail of distribution for the normalized sums of
multi-indexed independent non-negative random variables from (1.1). } \par

\vspace{4mm}

 \ Offered here results are generalizations of many ones obtained by S.Klesov  in an articles [14]-[17], see also [18], [26]-[27],
 where was obtained in particular the CLT for the sums of centered multiple variables. \par

\vspace{4mm}

 \ The multidimensional  case, i.e. when $ \ \vec{k} \in  Z_+^d, \ $ will be considered further. \par

\vspace{4mm}

 \ The paper is organized as follows. In the second section we describe and investigate the notion of
 the degenerate functions and approximation. In the next sections we obtain one of the main results:
 the  moment estimates for the multi-index sums with non-negative kernel in the two-dimensional case. \par
\ The so-called non-rectangular case is considered in the  $\ 4^{th} \ $ section. The fifth section is devoted
to the multivariate case. In the following section we obtain the exponential bounds for distribution of positive multiple sums.\par
 \ We show in the seventh section the upper bounds for these statistics. The last section contains as ordinary the concluding remarks. \par

\vspace{4mm}

 \ We reproduce here for readers convenience some results concerning the one-dimensional case, see an article  [30]. \par

\vspace{4mm}

 \ Denote for any r.v. - s. $ \ \eta_j, \ j = 1,2,\ldots  \ $ its $ \ k^{th} \ $ absolute moment by  $ \ m_j(k):  \  $

$$
m_j(k) := {\bf E} |\eta_j|^k, \ k \ge 1;
$$
 so that

$$
|\eta_j|_k = \left[ \ m_j(k) \right]^{1/k}.
$$

 \ We deduce applying the triangle inequality for the $ \ L(p, \Omega) \ $ norm   a very simple estimation

$$
|\sum_{j=1}^n \eta_j|_p \le \sum_{j=1}^n \left[ \ m_j(p) \ \right]^{1/p},
$$
and we do not suppose wherein that the r.v. $ \ \eta_j \ $ are non-negative and independent. \par

 \ In order to  describe a more fine estimations, we need to introduce some notations. \par

\vspace{4mm}

 \  Let us define the so-called {\it Bell's function} of two variables as follows.

$$
B(p,\beta) \stackrel{def}{=} e^{-\beta} \sum_{k=0}^{\infty} \frac{k^p \ \beta^k}{k!}, \ p \ge 2, \ \beta > 0,
$$
and put $  \ B(p) = B(p,1),  \ $  so that

$$
B(p) \stackrel{def}{=} e^{-1} \sum_{k=0}^{\infty} \frac{k^p}{k!}, \ p \ge 0.
$$

 \ The sequence of the numbers $ \ B(0) = 1, B(1),  B(2), B(3), B(4),\ldots  \  $ are named as Bell numbers; they appears in combinatorics, theory of probability etc., see [30], [25].

\vspace{4mm}

 \ Let the random variable (r.v.)  $ \ \tau = \tau[\beta], \ $ defined on certain probability space $ \ (\Omega, F,{\bf P}) \ $ with expectation $ \ {\bf E}, \ $
 has a Poisson distribution with parameter $ \  \beta, \ \beta > 0; \   $ write $ \ \Law(\tau) = \Law \tau[\beta] = \Poisson(\beta):   \  $

$$
{\bf P}(\tau = k) = e^{-\beta} \frac{\beta^k}{k!}, \ k = 0,1,2,\ldots,
$$

 \ It is worth to note that

$$
B(p,\beta) = {\bf E} (\tau[\beta])^p, \ p \ge 0.\eqno(1.2)
$$

\ In detail, let $ \ \eta_j, \ j = 1,2,\ldots \ $ be a sequence of non - negative independent random (r.v.); the case of centered or moreover symmetrical distributed r.v.
was considered in many works, see e.g. [9], [14]-[17], [26]-[27], [30] and so one. \par
 \  The following inequality holds true

$$
{\bf E} \left( \sum_{j=1}^n \eta_j   \right)^p \le B(p) \ \max \left\{  \ \sum_{j=1}^n {\bf E} \eta_j^p, \ \left( \ \sum_{j=1}^n {\bf E} \eta_j \ \right)^p \ \right\}, \ p \ge 2, \eqno(1.3)
$$
where the "constant" $ \ B(p) \ $ in (1.2) is the best possible, see [30]. \par

 \ One of the interest applications of these estimates in statistics, more precisely,  in the theory of $ \ U \ $ statistics  may be found in the articles [9], [18]. \par

 \ Another application. Let $ \  n = 1,2,3,\ldots; \  a, b = \const > 0; \ p \ge 2, \ \mu = \mu(a,b;p) := a^{p/(p-1)} \ b^{1/(1-p)} . \ $
 Define the following class of the sequences of an independent non - negative random variables

$$
Z(a,b) \stackrel{def}{=} \left\{ \eta_j, \ \eta_j \ge 0, \ \sum_{j=1}^n {\bf E}\eta_j = a; \ \sum_{j=1}^n {\bf E}\eta_j^p = b \ \right\}. \eqno(1.4)
$$

 \ G.Schechtman in proved that

$$
\sup_{ \ n = 1,2,\ldots; \  \{\eta_j \} \in Z(a,b) \  } {\bf  E}  \left( \ \sum_{j=1}^n \eta_j \  \right)^p = \left(  \frac{b}{a} \right)^{p/(p-1)} \ B(p, \mu(a,b;p)). \eqno(1.5)
$$

\vspace{4mm}

\ The introduced above Bell's function allows in turn a simple estimation, see [30], which may be used by practical using.
Indeed, let us introduce the following  auxiliary function

$$
g_{\beta}(p) \stackrel{def}{=} \frac{p}{e} \cdot
 \inf_{\lambda > 0} \left[\lambda^{-1} \ \exp \left( \ \beta \left(e^{\lambda} - 1  \right) \ \right) \ \right]^{1/p}, \ \beta, p  > 0.  \eqno(1.6)
$$

 \ It is proved in particular in  [30] that

$$
B^{1/p}(p, \beta) \le  g_{\beta}(p), \ p,\beta > 0. \eqno(1.7)
$$

\vspace{4mm}

\ Let us introduce also the following function

$$
h_0(p, \beta) \stackrel{def}{=} \sup_{k=1,2,\ldots}  e^{-\beta}  \left\{ \ \frac{k^p \ \beta^k}{k!} \ \right\}; \eqno(1.8)
$$
therefore

$$
B(p,\beta) \ge h_0(p, \beta), \ p, \ \beta > 0. \eqno(1.9)
$$

\vspace{4mm}

 \ The last estimate may be simplified in turn as follows. We will apply the following version of the famous Stirling's formula

$$
 k! \le \zeta(k), \ \ k = 1,2,\ldots,
$$
where

$$
\zeta(k) \stackrel{def}{=} \sqrt{2 \pi k} \ \left( \frac{k}{e}  \right)^k \ e^{1/(12k)}, \ k = 1,2,\ldots \eqno(1.10)
$$

 \ Define a new function

$$
h(p, \beta) \stackrel{def}{=} \sup_{x \in(1, \infty)} \
 \left\{ e^{ 1/(6 p x )} \cdot \left[ \ \frac{e^{x - \beta} \ x^{p - x - 1/2}}{ \sqrt{2 \ \pi} \ x^x}  \ \right]^{1/p} \right\}. \eqno(1.11)
$$

 \ We obtained in [30] really the following lower simple estimate for the Bell's function \par

\vspace{4mm}

$$
B^{1/p} (p, \beta) \ge h_0(p, \ \beta), \  B^{1/p} (p, \beta) \ge h(p, \ \beta), \ p, \ \beta > 0. \eqno(1.12)
$$

\vspace{4mm}

 \ These estimates may be in turn simplified as follows. Assume that $ \ p \ge 2 \beta, \ \beta > 0, \ p \ge 1; \ $ then

$$
B^{1/p}(p,\beta) \le \frac{p/e}{\ln (p/\beta) - \ln \ln (p/\beta)} \cdot \exp \left\{ \frac{1}{\ln (p/\beta)}  - \frac{1}{p/\beta}  \right\}. \eqno(1.13)
$$

\vspace{4mm}

 \ For example,

$$
B^{1/p}(p) \le \frac{p/e}{\ln p - \ln \ln p} \cdot \exp \left\{ 1/\ln p - 1/p   \right\}, \ p \ge 2. \eqno(1.14)
$$

\vspace{4mm}

\ The estimate (1.13) may be simplified as follows

$$
B^{1/p}(p) \le \frac{p}{e \ \ln (p/\beta)} \cdot \left[ 1 + C_1(\beta) \cdot \ \frac{\ln \ln (p/\beta)}{\ln (p/\beta)} \  \right], \eqno(1.15)
$$
where $  \ C_1(\beta) = \const \in (0,\infty), $ and we recall that $  \ p \ge 1, \ p \ge 2\beta. \ $ \par

\vspace{4mm}

 \ For example,

$$
B^{1/p}(p) \le \frac{p/e}{\ln p - \ln \ln p} \cdot \exp \left\{ 1/\ln p - 1/p   \right\}, \ p \ge 2. \eqno(1.16)
$$

\vspace{4mm}

\ The lower estimate for Bell's function has a form

 $$
 \ B^{1/p}(p,\beta)  \ge \
 $$

$$
 \beta^{1/\ln(pe/\beta)} \cdot \frac{p}{\ln (pe/\beta)} \cdot \left\{ \exp \left[ \ \frac{\ln p - \ln(pe)/\beta}{\ln[(pe)/\beta]} \ \right] \  \right\}^{-1},
$$

$$
 \ p,\beta > 0, p/\beta \ge 2.
$$

 \ It may be simplified as follows

$$
B^{1/p}(p) \ge \frac{p}{e \ \ln (p/\beta)} \cdot \left[ 1 - C_2(\beta) \cdot \ \frac{\ln \ln (p/\beta)}{\ln (p/\beta)} \  \right], \eqno(1.17)
$$
where
$  \ C_2(\beta) = \const \in (0,\infty), $ and we recall that $ \ p \ge 2\beta.\ $

\vspace{4mm}

 \ We suppose hereafter that both the variables $ \ p \ $ and $ \ \beta \ $ are independent but such that

$$
p \ge 1, \ \beta > 0, \ p/\beta  \le 2. \eqno(1.18)
$$

 \ It is known [30] that in this case  there exist   {\it two absolute positive constructive finite constants} $ \ C_3, \ C_4, \ C_3 \le C_4 \ $ such that

$$
C_3 \ \beta \le B^{1/p} (p, \beta) \le C_4 \ \beta. \eqno(1.19)
$$

\vspace{4mm}

 \begin{center}

 {\it To summarize.}

\end{center}

 \ Define the following infinite dimensional random vector

$$
\eta = \vec{\eta} = \{ \eta_1, \eta_2, \eta_3, \ldots    \}. \eqno(1.20)
$$
 \ Recall that here the r.v. $ \ \{ \eta_j \} \ $ are non-negative and independent. One can suppose also that $ \ m_j(p) = |\eta_j|_p < \infty \ $ for some value $ \ p \in [2, \infty). \ $ \par
 \ Put  also

 $$
 Z_{p,n} = Z_{p,n}[\eta]  := n^{-1} \sum_{j=1}^n \left[ \ m_j(p) \ \right]^{1/p},
 $$

$$
V_{p,n} = V_{p,n}[\eta] := n^{-1} \ \cdot B^{1/p}(p) \cdot \ \max \left\{ \left( \sum_{j=1}^n m_j(p) \right)^{1/p}, \ \left( \ \sum_{j=1}^n m_j(1)  \right) \  \right\},
$$
and ultimately

$$
\Theta_{p,n} = \Theta_{p,n}[\eta]  :=\min \{  Z_{p,n}[\eta],  V_{p,n}[\eta]  \}, \eqno(1.21)
$$

$$
\Theta_{p} = \Theta_p[\eta] :=  \sup_{n = 1,2,3,\ldots}   \Theta_{p,n}, \eqno(1.22)
$$
with described above correspondent upper estimations for the function $ \ B^{1/p}(p).  \  $\par

\vspace{4mm}

{\bf Proposition 1.1.} We deduce under formulated conditions

$$
 \left| n^{-1} \ \sum_{j=1}^n \eta_j  \right|_p \le \Theta_{p,n}[\eta]. \eqno(1.23)
$$
 \ As a slight consequence

$$
\sup_{n=1,2,3,\ldots} \left| n^{-1} \ \sum_{j=1}^n \eta_j  \right|_p \le \Theta_{p}[\eta]. \eqno(1.23a)
$$

\vspace{4mm}

{\bf Remark 1.1.} Of course, one can use in the practice in the estimate (1.21) and in the next  relations instead the variable $ \  B^{1/p}(p)  \  $ its
upper bounds from the inequality (1.16). \par

\vspace{4mm}

 \section{Degenerate functions and approximation. }

\vspace{4mm}

 \ {\bf Definition 2.1,} see [26], [27], [30].  The measurable  function $ \  f:  X \otimes Y \to R  \ $ is said to be {\it  degenerate, } if it has a form

$$
f(x,y) = \sum_{k(1) = 1}^M \sum_{k(2)=1}^M \lambda_{k(1), k(2)} \ g_{k(1)}(x) \ h_{k(2)}(y),  \eqno(2.1)
$$
where $ \  \lambda_{i,j} = \const \in R, \  M =  \const = 1,2, \ldots, \infty. \ $ \par

 \ One can distinguish two cases in the relation (2.1): ordinary, or equally {\it finite degenerate function,} if in (2.1) $ \ M < \infty, \ $  and
 infinite degenerate function otherwise. \par

  \ The degenerate functions (and kernels) of the form (2.1) are used, e.g., in the approximation theory, in the theory of random processes and fields,
in the theory of integral equations, in the game theory etc. \par

 \  A particular application of this notion may be found in the  authors articles [26], [27], [30]. \par

\vspace{4mm}

 \ Denotation: $ \ M = M[f] \stackrel{def}{=}  \deg(f); \ $ of course, as a capacity of the value $ \ M \ $  one can understood its {\it constant}
minimal value. \par

 \ Two examples. The equality (2.1) holds true if the function $ \  f(\cdot, \cdot) \  $  is trigonometrical or algebraical polynomial.  \par
\ More complicated example: let $  \   X \ $ be compact metrisable space  equipped with the non-trivial probability Borelian measure $ \ \mu. \ $
This  imply that an arbitrary non-empty open set has a positive measure. \par

 \ Let also $ \  f(x,y), \  x,y \in X \ $ be continuous numerical valued  non-negative definite function. One can write the famous Karunen-Loev's
decomposition

$$
f(x,y) = \sum_{k=1}^M  \lambda_k  \ \phi_k(x) \ \phi_k(y),\eqno(2.2)
$$
where $ \  \lambda_k, \ \phi_k(x) \ $ are  correspondingly eigen values and eigen orthonormal
function for the function (kernel)  $ \ f(\cdot, \cdot): \ $

$$
\lambda_k \ \phi_k(x)  = \int_X f(x,y) \phi_k(y) \ \mu(dy).
$$
 \ We assume without loss of generality

$$
\lambda_1 \ge \lambda_2 \ge \ldots \lambda_k \ge \ldots \ge 0.
$$

\vspace{4mm}

 \ {\it  It will be presumed in this report, i.e. when  the function $ \ f = f(x,y) \ $  is non negative,
 in addition to the expression (2.1) that all the functions   \ } $ \  \{  g_i \}, \ \{ h_j \}  \  $
{\it are also non-negative:}

$$
\forall x \in X \ g_i(x) \ge 0, \ \forall y \in Y \ h_j(y) \ge 0. \eqno(2.3)
$$

\  Further, let  $ \  B_1, \ B_2, \  B_3, \ \ldots, B_M \  $  be some rearrangement invariant (r.i.) spaces builded correspondingly over the
spaces $ \ X,Y; \ Z,W, \ldots, \ $
for instance, $ \ B_1 = L_p(X), \  B_2 = L_q(Y), 1 \le p,q \le \infty.  \  $ If $ \  f(\cdot) \in B_1 \otimes B_2, \  $ we suppose also
in (2.1) $ \ g_i \in B_1, \ h_j \in B_2; \ $ and if in addition in (2.1) $ \  M = \infty,  \  $ we suppose that the series in (2.1)
converges in the norm $  \  B_1 \otimes B_2 \   $

$$
\lim_{m \to \infty} || \   f(\cdot) -   \sum_{i = 1}^m \sum_{j=1}^m \lambda_{i,j} \ g_i(\cdot) \ h_j(\cdot)  \  ||B_1 \otimes B_2 = 0.
\eqno(2.3b)
$$

\  The condition (2.3b) is satisfied if for example $ \ ||g_i||B_1  = ||h_j||B_2 = 1 \ $ and

$$
 \sum \sum_{i,j = 1}^M |\lambda_{i,j}| < \infty, \eqno(2.4)
$$
or more generally when

$$
\sum \sum_{i,j = 1}^M |\lambda_{i,j}| \cdot ||g_i||B_1 \cdot ||h_j||B_2 < \infty. \eqno(2.4a)
$$

 \ The function of the form (2.1) with $ \    M = M[f] = \deg (f) < \infty \ $ is named {\it degenerate},
notation $ \  f \in D[M]; \ $ we put also $  \  D := \cup_{M < \infty} D[M]. \ $ Obviously,

$$
B_1 \otimes B_2 = D[\infty].
$$

 \ Define also for each  {\it non-negative} such a function $ \   f \in D \  $ the following  quasi-norm, also non-negative:

$$
||f|| D^+(B_1, B_2) \stackrel{def}{=} \inf \left\{ \  \sum \sum_{i,j = 1,2,\ldots,M[f] }  |\lambda_{i,j}| \ ||g_i||B_1 \ ||h_j||B_2 \  \right\}, \eqno(2.5)
$$
where all the arrays $ \ \{ \lambda_{i,j} \} , \ \{ g_i\}, \ \{h_j\}  \ $ are taking from the  representation  2.1; and in addition,

$$
 g_i(x) \ge 0, \ h_j(y) \ge 0; \  x \in X, \ y \in Y.
$$

 \ We will write for brevity $ \ ||f||D_p   := $

$$
||f|| D^+(L_p(X), L_p(Y)) =
\inf \left\{ \  \sum \sum_{i,j = 1,2,\ldots,M[f] }  |\lambda_{i,j}| \ |g_i|_p \ |h_j|_p \  \right\}, \eqno(2.5a)
$$
where all the arrays $ \ \{ \lambda_{i,j} \} , \ \{ g_i\}, \ \{h_j\}  \ $ are taking from the  representation  2.1, of course with non - negative functions $ \ g_i, \ h_j. \ $  \par

\vspace{4mm}

\ Further, let the function  $ \  f \in B_1 \otimes B_2  \ $ be given. The error of a degenerate approximation of the non-negative function $ \ f: X \otimes Y \to R \ $
by the degenerate ones of the degree $ \ M, \ $  also with non-negative summands, will be introduced as follows

$$
Q^+_M[f](B_1 \otimes B_2) \stackrel{def}{=} \inf_{\tilde{f} \in D[M]} ||f - \tilde{f}||B_1 \otimes B_2 =
\min_{\tilde{f} \in D^+[M]} ||f - \tilde{f}||B_1 \otimes B_2. \eqno(2.6)
$$
 \  Obviously, $ \  \lim Q^+_M[f] (B_1 \otimes B_2) = 0, \ M \to \infty. \  $\par

 \ For brevity:

$$
Q^+_M[f]_p \stackrel{def}{=} Q^+_M[f](L_p(X) \otimes L_p(Y)). \eqno(2.6a)
$$

\vspace{4mm}

 \  The function $ \ \tilde{f} \ $ which realized the minimum in (2.6), obviously, non-negative, not necessary to be unique,
will be  denoted by $ \  Z^+_M[f](B_1 \otimes B_2):  \  $

$$
 Z^+_M[f](B_1 \otimes B_2):= \argmin_{\tilde{f} \in D^+[M]} ||f - \tilde{f}||B_1 \otimes B_2, \eqno(2.7)
$$
so that

$$
Q^+_M[f](B_1 \otimes B_2) = ||f - Z^+_M[f]||(B_1 \otimes B_2).  \eqno(2.8)
$$

 \ For brevity:

$$
Z^+_M[f]_p := Z^+ _M [f](L_p(X) \otimes L_p(Y)). \eqno(2.9)
$$

 \ Let for instance again $ \  f(x,y), \  x,y \in X \ $ be continuous numerical valued  {\it non-negative definite} function, non necessary to be non - negative,
see (2.3) and (2.3a). It is easily to calculate

$$
Q_M[f] (L_2(X) \otimes L_2(X)) = \sum_{k=M + 1}^{\infty} \lambda_k.
$$

\vspace{4mm}

 \section{Moment estimates for multi-index sums.}

\vspace{4mm}

\begin{center}

 {\bf  Two-dimensional case. } \par

\end{center}

\vspace{4mm}

{\bf \  A trivial estimate.  }  \par

\vspace{4mm}

 \ The following simple estimate based only on the triangle inequality, may be interpreted as trivial:

$$
|S_L|(L_p(X) \otimes L_p(Y)) \le \ |f|(L_p(X) \otimes L_p(Y)) = |f|_p, \eqno (3.1)
$$
even without an assumption of the non-negativity of the function $ \ f \ $ and the independents of the r.v $ \  g_k(\xi(i)), \  h_l(\eta(j)). \ $ \par

\vspace{4mm}

{\it  Hereafter \   } $ \ p \ge 2. \ $ \par

\vspace{4mm}

\ {\bf   The two-dimensional degenerate case.  } \par

\vspace{4mm}

 \ In this subsection  the non-negative kernel-function $ \  f = f(x,y) \  $ will be presumed to be degenerate with minimal constant possible
degree $ \ M = M[f] = 1, \ $ on the other words, {\it factorizable}:

$$
f(x,y) = g(x) \cdot h(y), \ x \in X, \ y \in Y, \eqno(3.2)
$$
of course, with non-negative factors $ \ g,h. \ $ \par

 \ Further, we suppose that the set $ \ L \ $ is integer constant rectangle:

$$
L = [1,2,\ldots, n(1)] \otimes [1,2, \ldots, n(2)], \ n(1), n(2) \ge 1.
$$

 \ Let us consider the correspondent double sum $  \  S_L[f] =  S^{(2)}_L[f] :=   $

$$
 |L|^{-1} \sum  \sum_{i, j \in L}  g(\xi_i) \ h(\eta_j), \ \eqno(3.3),
$$
where
$$
 n =  \vec{n} = (n(1), n(2)) \in L, \ n_1, n_2 \ge 1. \eqno(3.3)
$$

  \ We have denoting

$$
\vec{g} = \{ g(\xi(i)) \}, \ i = 1,2, \ldots, n(1); \ \vec{h}  = \{ h(\eta(j)) \}, \ j = 1,2,\ldots, n(2):
$$

$$
 S_L[f] = \left[ n(1)^{-1} \sum_{i=1}^{n(1)} g(\xi(i)) \right] \cdot \left[ n(2)^{-1} \sum_{j=1}^{n(2)} h(\eta(j)) \right].
$$
 \ Since both the factors in the right-hand size of the last inequality are independent, we deduce
 applying the one-dimensional estimates (1.23), (1.23a):

$$
|S_L|_p \le \Theta_{p, n(1)}[\vec{g}] \cdot \Theta_{p, n(2)}[\vec{h}], \eqno(3.4)
$$
and hence

$$
\sup_{L: |L| \ge 1} |S_L|_p \le \Theta_{p}[\vec{g}] \cdot \Theta_{p}[\vec{h}], \eqno(3.4a)
$$

\vspace{4mm}

{\bf  Estimation for an arbitrary degenerate kernel.  } \par

\vspace{4mm}

 \ In this subsection the function $ \ f(\cdot,\cdot) \ $ is non-negative and degenerate, as well as all the functions $ \ g_{k(1)}(\cdot), \ h_{k(2)}(\cdot): \ $

$$
f(x,y) = \sum \sum_{k(1), k(2) = 1}^M \ \lambda_{k(1), k(2)} \ g_{k(1)}(x) \ h_{k(2)}(y), \eqno(3.5)
$$
where $ \ 1 \le M \le \infty, \ $

$$
g_{k(1)}(\cdot) \in L_p(X), \ h_{k(2)} (\cdot) \in L_p(Y),
$$
and as before $ \  g_{k(1)}(x) \ge 0, \  h_{k(2)}(y) \ge 0.  \  $ \par

\vspace{4mm}

 \ Denote also by $ \ R = R(f)  \  $ the {\it set } of all such the functions $ \  \{g\} = \vec{g}   \  $ and  $ \ \{h\} = \vec{h} \ $  as well as the sequences of coefficients
  $  \  \{\lambda\} = \{ \ \lambda_{k(1), k(2)} \ \}  \  $ from the representation (3.5):

$$
R[f] := \{  \{\lambda\}, \  \{g\} = \vec{g}, \ \{h\} = \vec{h} \}:
$$

$$
 f(x,y) = \sum \sum_{k(1), k(2) = 1}^M \ \lambda_{k(1), k(2)} \ g_{k(1)}(x) \ h_{k(2)}(y). \eqno(3.5a)
$$

\vspace{4mm}

 \ {\it We must impose on the series (3.5) in the case  when $ \ M = \infty \ $  the condition of its convergence in the norm of the space} $ \  L_p(X) \otimes L_p(Y). \  $\par

\vspace{4mm}

 \ Let us investigate the introduced before  statistics

$$
S^{(\lambda)}_{L} = S^{(\lambda)}_L[f]  := | L|^{-1} \ \sum \sum_{i,j \in L} f(\xi(i), \eta(j)) =
$$

$$
| L|^{-1} \ \sum \sum_{i,j \in L}  \left[ \ \sum \sum_{k(1), k(2) = 1}^M \lambda_{k(1), k(2)} \ g_{k(1)}(\xi(i)) \ h_{k(2)}(\eta(j)) \ \right] =
$$

$$
|L|^{-1} \ \sum \sum_{k(1), k(2) = 1}^M \lambda_{k(1), k(2)} \ \left[ \ \sum \sum_{i,j  \in L} g_{k(1)}(\xi_i) \ h_{k(2)}(\eta_j) \ \right]=
$$

$$
 \sum \sum_{k(1), k(2) = 1}^M \lambda_{k(1), k(2)} \cdot \left[ \ (n(1))^{-1} \sum_{i=1}^{n(1)} g_{k(1)}(\xi(i)) \ \right]  \times
$$

$$
  \left[ (n(2))^{-1} \sum_{j=1}^{n(2)} h_{k(2)}(\eta(j)) \right]. \eqno(3.6)
$$

 \  We have using the triangle inequality and the estimate (3.4)

$$
\left| S^{(\lambda)}_L[f] \right|_p \le \sum \sum_{k(1), k(2) = 1}^M | \ \lambda_{k(1), \ k(2)} \ | \ \cdot
 \Theta_{p, n(1)} \left[ \vec{g}_{k(1)} \ \right] \cdot \Theta_{p, n(2)} \left[ \ \vec{h}_{k(2)} \  \right]. \eqno(3.7)
$$

 \ This estimate remains true in the case when $ \ M = \infty, \ $ if of course the right - hand side of (3.7) is finite; in the opposite case it is nothing to  make.  \par

 \vspace{4mm}

  \ To summarize, we introduce a new weight norm on the (numerical) array $ \ \vec{\lambda} = \{ \ \lambda_{k(1), \ k(2)} \}, $ more exactly, the sequence of the norms

$$
f \in D(M) \ \Rightarrow |||f||| \Theta_p = |||\vec{\lambda}|||\Theta_p =
$$

$$
  |||\vec{\lambda}|||\Theta_p^{(2)} = |||\vec{\lambda}|||\Theta(p;n(1), n(2), \ \{g\}, \ \{h\}) \stackrel{def}{=}
$$

$$
\sum \sum_{k(1), k(2) = 1}^M | \ \lambda_{k(1), \ k(2)} \ | \ \cdot
 \  \Theta_{p, n(1)} \left[ \ \vec{g}_{k(1)} \ \right] \cdot \Theta_{p, n(2)} [ \ \vec{h}_{k(2)} \ ].\eqno(3.8)
$$

\vspace{4mm}

 {\bf Proposition 3.1.} If $ \ f \in D(M), \ $ then

$$
\left| S^{(\lambda)}_L[f] \right|_p \le |||f||| \Theta_p  = |||\vec{\lambda}|||\Theta_p =
$$

$$
 |||\vec{\lambda}|||\Theta(p;n(1), n(2),  \ \{g\}, \ \{h\}), \eqno(3.9)
$$
and as a consequence

$$
 \sup_{ L: |L| \ge 1} \left| S^{(\lambda)}_L[f] \right|_p \le \sup_{n(1), n(2)}  |||f||| \Theta_p =\sup_{n(1), n(2)} |||\vec{\lambda}|||\Theta_p\{f\} =
$$

$$
  \sup_{n(1), n(2)} |||\vec{\lambda}|||\Theta(p;n(1), n(2), \{g\}, \{h\}). \eqno(3.9a)
$$

\vspace{4mm}

{\bf Main result.  Degenerate approximation approach. } \par

\vspace{4mm}

 \ {\bf Theorem 3.1.}  Let $ \  f = f(x,y)  \  $ be arbitrary function from the space $ \   L_p(X) \otimes L_p(Y), \ p \ge 2. \  $ Then
$ \ |S_L[f]|_p \le W_L[f](p), $ where

\vspace{4mm}

$$
\ W_L[f](p) =  W_L[f; \{\lambda\}, \ \{g\}, \ \{h\} \ ](p) \stackrel{def}{=}
$$

$$
\inf_{M \ge 1} \left[ \ |||Z_M[f]|||\Theta_p  +  \ Q^+_M[f]_p \  \right], \eqno(3.10)
$$
where in turn the vector triple $ \ \{\lambda\}, \ \{g\}, \ \{h\} \ $ is taken from the representation (3.5):  $ \left[ \{\lambda\}, \ \{g\}, \ \{h\} \ \right] \in R[f]. \ $ \par

\vspace{4mm}

 \ As a slight consequence: $ \ \sup_{L: |L| \ge 1} |S_L[f]|_p \le W[f](p), $ where

$$
 \ W[f](p) \stackrel{def}{=} \sup_{L: |L| \ge 1} W_L[f](p) =
$$

$$
 \sup_{L: |L| \ge 1} \inf_{M \ge 1} \left[ \ |||Z_M[f]|||\Theta_p + \ Q^+_M[f]_p \  \right], \eqno(3.10a)
$$
and the next consequence

 $$
 \sup_{L: |L| \ge 1} |S_L[f]|_p \le \inf_{ \{\lambda\}, \ \{g\}, \ \{h\} \in R[f] } W_L[f; \{\lambda\}, \ \{g\}, \ \{h\} \ ](p). \eqno(3.10b)
 $$

\vspace{4mm}

{\bf Proof } is very simple, on the basis of previous results of this section. Namely, let $ \ L \ $ be an arbitrary non-empty set. Consider a
splitting

$$
f = Z^+_M[f] + ( f - Z^+_M[f] ) =: \Sigma_1 + \Sigma_2.
$$
 \ We have

$$
|\Sigma_1|_p   \le |||Z_M[f]|||\Theta_p.
$$

 \ The member $  \ |\Sigma_2|_p $ may be estimated by virtue of inequality (3.1):

$$
|\Sigma_2|_p \le |L| \ |f - Z^+_M[f]|_p =  \ Q^+_M[f]_p.
$$
 \ It remains to apply the triangle inequality and minimization over $ \ M. \ $ \par

\vspace{4mm}

 \ {\bf Example 3.1.}  We deduce from (3.10) as a particular case

$$
 \sup_{ L: |L| \ge 1} \left| S^{(\lambda)}_L[f] \right|_p \le  \  ||f||D^+_p, \eqno(3.11)
$$
 if of course the right-hand side  of (3.11) is finite for some value $ \ p, \ p \ge 2. \ $ \par
\ Recall that in this section $ \ d = 2. \ $ \par

 \vspace{4mm}

 \section{ Non-rectangular case.   } \par

 \vspace{4mm}

 \ We denote  by $ \  \pi^+(L) \  $ the set of all rectangular's  which are{\it\ circumscribed\ }about
 the set $ \ L: \ \pi^+(L) = \{ L^+  \}, $  where

$$
L^+ = \{ [n(1)^{+}, n(1)^{++}] \otimes [n(2)^{+}, n(2)^{++}]  \}: \ L^+ \supset L, \eqno(4.1)
$$
and

$$
1 \le n(1)^{+}  \le n(1)^{++} < \infty, \   n(1)^{+},  n(1)^{++} \in Z_+,
$$

$$
1 \le n(2)^{+}  \le n(2)^{++} < \infty, \   n(2)^{+},  n(2)^{++} \in Z_+.
$$

\vspace{4mm}

\ {\bf Proposition 4.1.}

$$
| S_L|_p \le \inf_{L^+: \ L \subset L^+} \  \left\{ \ \frac{|L^+|}{|L|} \cdot  W_{L^+}[f; \{\lambda\}, \ \{g\}, \ \{h\} \ ](p) \ \right\}. \eqno(4.2)
$$

\vspace{4mm}

 \ {\bf  Proof} is very simple. We have

$$
|L| \ S_L = \sum \sum_{i,j \in L} f(\xi(i), \ \eta(j)),
$$
therefore

$$
\left| \ |L| \ S_L \ \right|_p  =  \left| \sum \sum_{i,j \in L} f(\xi(i), \ \eta(j)) \ \right|_p \le
$$

$$
\left| \sum \sum_{i,j \in L^+} f(\xi(i), \ \eta(j)) \ \right|_p \le |L^+| \cdot W_{L^+}[f: \{\lambda\}, \ \{g\}, \ \{h\} \ ](p),
$$
by virtue of theorem 3.1. \par

\vspace{4mm}

 \section{ Moment estimates for multi-index sums.}

 \vspace{4mm}

\begin{center}

{\bf  Multidimensional generalization.}\par

\end{center}

\vspace{4mm}

 \ Let now $ \  (X_m, \ B_m, \mu_m), \ m = 1,2,\ldots, d, \ d \ge 3  \ $ be {\it a family of probability spaces:}  $ \  \mu_m(X_m)  = 1; \  $
$ \ X := \otimes_{m=1}^d X_m; \ \xi(m) \ $ be independent random variables having the distribution correspondingly
$ \ \mu_m: \ {\bf P}(\xi(m) \in A_m) = \mu_m(A_m), \ A_m \in B_m; \ $
$  \  \xi_i(m), \ i = 1,2,  \ldots, n(m); \ n(m) = 1,2,  \ldots, \ n(m) < \infty  \  $ be independent  copies of $ \ \xi(m) \ $ and also independent
on the other vectors  $ \  \xi_i(s), s \ne m, \ $ so that all the random variables $ \ \{ \xi_i(m)  \} \ $ are common independent. \par

 \ Another notations, conditions, restrictions and definitions. $ \  L  \subset Z_+^d, \ |L| = \card(L) > 1;  \ j = \vec{j} \in L; \  $

$$
 k = \vec{k} = (k(1), k(2), \ldots, k(d)) \in Z_+^d; \  N(\vec{k}) := \max_{j = 1,2, \ldots,d} k(j); \eqno(5.0)
$$

$ \vec{\xi} := \{\xi(1), \xi(2), \ldots, \xi(n(m)) \}; \   \vec{\xi}_i := \{\xi_i(1), \xi_i(2), \ldots, \xi_i(n(m)) \};  $
$ \ X := \otimes_{i=1}^d X_i, \ f:X \to R \  $ be measurable {\it non-negative} function, i.e. such that $ \  f(\vec{\xi}) \ge 0; \ $

$$
S_L[f] := |L|^{-1} \sum_{k \in L} f\left(\vec{\xi}_k \right). \eqno(5.1)
$$

 \ The following simple estimate is named as before trivial:

$$
|S_L[f]|_p  \le  |f|L_p. \eqno (3.0a)
$$

\vspace{4mm}

{\it  Recall that  by-still hereafter \   } $ \ p \ge 2. \ $ \par

\vspace{4mm}

 \ By definition, as above, the function $ \  f: X \to R  \ $ is said to be degenerate, iff it has the form

$$
f(\vec{x}) = \sum_{\vec{k} \in Z_+^d, \ N(\vec{k}) \le M} \lambda(\vec{k}) \ \prod_{s=1}^d g^{(s)}_{k(s)}(x(s)), \eqno(5.2)
$$
for some integer {\it constant} value $ \ M, \ $  finite or not, where all the functions $ \   g^{(s)}_k(\cdot) \  $ are in turn non-negative:
$ \  g^{(s)}_k(\xi(k)) \ge 0. \ $  Denotation: $ \ M = \deg[f]. $ \par

\vspace{4mm}

 \ Define also as in the two-dimensional case for each such a function $ \   f \in D^+ \  $ the following non-negative quasi-norm

$$
||f|| D^+_p \stackrel{def}{=} \inf \left\{ \ \sum_{\vec{k} \in Z^d_+, \ N(\vec{k}) \le M[f] }
|\lambda(\vec{k})|  \cdot \prod_{s=1}^d  |g^{(s)}_{k(s)}(\xi(s))|_p \  \right\}, \eqno(5.3)
$$
where all the arrays $ \ \{ \lambda( \vec{k}) \} , \ \{ g_j \},  \ $ are taking from the  representation 5.2, in particular, all the summands are non-negative. \par

 \ The last assertion allows a simple estimate: $ \  ||f||D^+_p \le || f ||D^{+o}_p, \  $ where

$$
||f|| D^{+o}_p \stackrel{def}{=}  \ \sum_{\vec{k} \in Z^d_+, \ N(\vec{k}) \le M[f] }
|\lambda(\vec{k})|  \cdot \prod_{s=1}^d  |g^{(s)}_{k(s)}(\xi(s))|_p, \eqno(5.3a)
$$
and if we denote

$$
G(p) :=  \prod_{j=1}^d  |g^{(j)}_{k_j}(\xi_j)|_p, \ p \ge 1; \ || \lambda ||_1 := \sum_{\vec{k} \in Z^d_+}  |\lambda(\vec{k})|,
$$
then

$$
||f||D^+_p \le ||f||D^{+o}_p \le G(p) \cdot ||\lambda||_1. \eqno(5.3b)
$$

\vspace{4mm}

\ Further, let the non-negative  function  $ \  f \in B_1 \otimes B_2  \otimes \ldots \otimes B_d \ $ be given.
Here $ \ B_r, \ r = 1,2,\ldots,d \ $ are some Banach functional rearrangement invariant spaces builded correspondingly over the sets $ \ X_m. \ $ \par
 The error of a degenerate approximation
of the function $ \ f \ $ by the degenerate and non-negative ones of the degree $ \ M \ $ will be introduced as before

$$
Q^+_M[f](B_1 \otimes B_2 \otimes \ldots B_d) \stackrel{def}{=} \inf_{\tilde{f} \in D^+[M]} ||f - \tilde{f}||B_1 \otimes B_2 \otimes \ldots B_d =
$$

$$
\min_{\tilde{f} \in D^+M]} ||f - \tilde{f}||B_1 \otimes B_2 \otimes \ldots \otimes B_d. \eqno(5.4)
$$
  \ Obviously, $ \  \lim Q_M[f] (B_1 \otimes B_2 \otimes \ldots \otimes B_d) = 0, \ M \to \infty. \  $\par

 \ For brevity:

$$
Q_M[f]_p \stackrel{def}{=} Q_M[f](L_p(X_1) \otimes L_p(X_2) \otimes \ldots \otimes L_p(X_d)). \eqno(5.5)
$$

\vspace{4mm}

 \  The function $ \ \tilde{f} \ $ which realized the minimum in (5.4), not necessary to be unique,
will be denoted by $ \  Z_M[f](B_1 \otimes B_2 \otimes \ldots \otimes B_d) = Z^+_M[f](B_1 \otimes B_2 \otimes \ldots \otimes B_d):  \  $

$$
 Z^+_M[f](B_1 \otimes B_2 \otimes \ldots \otimes B_d):=
\argmin_{\tilde{f} \in D^+[M]} ||f - \tilde{f}||B_1 \otimes B_2 \otimes \ldots \otimes B_d, \eqno(5.6)
$$
so that

$$
Q^+_M[f](B_1 \otimes B_2 \otimes \ldots \otimes B_d) = ||f - Z^+_M[f]||B_1 \otimes B_2 \otimes \ldots \otimes B_d.  \eqno(5.7)
$$

 \ For brevity:

$$
Z^+_M[f]= Z_M[f]_p := Z_M [f](L_p(X_1) \otimes L_p(X_2) \otimes \ldots \otimes L_p(X_d) ). \eqno(5.8)
$$

\vspace{4mm}

 \ Denote as in the third section for $ \ f \in D(M) \ $ in the multivariate d-dimensional case

$$
|||f|||\Theta_p = |||f|||\Theta_{p,L} \stackrel{def}{=}
$$

$$
\sum_{N(\vec{k}) \le M} \left|\lambda_{\vec{k}} \right|  \ \prod_{s=1}^d \Theta_{p, n(s)} \left[g^{(s)}_{k(s)} \right], \eqno(5.9a)
$$

$$
W_L[f](p) = W_L^{(d)}[f](p)  \stackrel{def}{=}  \inf_M \left[ \  |||Z^+_M[f] + Q_M^+[f]_p  \ \right], \eqno(5.9b)
$$

$$
  W[f]^{(d)}(p) \stackrel{def}{=} \sup_{L: |L| \ge 1} W_L^{(d)}[f](p). \eqno(5.9c)
$$

\ We deduce analogously to the third section \par

\vspace{4mm}

 {\bf Proposition 5.1.} If $ \ f \in D(M), \ $ then

$$
\left| S^{(\lambda)}_L[f] \right|_p \le |||f||| \Theta_{p,L}, \eqno(5.10)
$$
and of course

$$
\sup_{L: |L| \ge 1} \left| S^{(\lambda)}_L[f] \right|_p \le  \sup_{L: |L| \ge 1}  |||f||| \Theta_{p,L}, \eqno(5.10a)
$$

\vspace{4mm}

 \ {\bf Theorem 5.1.}  Let $ \  f = f(x) = f(\vec{x}), \ x \in X  \  $ be arbitrary non-negative function from the space
$ \   L_p(X_1) \otimes L_p(X_2) \otimes \ldots \otimes L_p(X_d), \ p \ge 2. \  $ Then

$$
|S_L[f]|_p \le  W_L^{(d)}[f](p),\eqno(5.11)
$$

$$
 \sup_{L: |L| \ge 1} |S_L[f]|_p \le   \sup_{L: |L| \ge 1} W_L^{(d)}[f](p) =  W^{(d)}[f](p) . \eqno(5.11a)
$$

\vspace{4mm}

 \ {\bf Example 5.1.}  We deduce  alike the example 3.1 as a particular case

$$
 \sup_{ L: |L| \ge 1} \left| S_L[f] \right|_p \le  (2/3)^d \cdot \left[ \ \frac{p}{e \cdot \ln p} \ \right]^d \cdot ||f||D_p,  \eqno(5.12)
$$
 if of course the right-hand side  of (5.9a) is finite for some value $ \ p, \ p \ge 2. \ $ \par

\vspace{4mm}

 \ {\bf Remark 5.1.} Notice that the last estimates (5.10), (5.11), and (5.12) are essentially non-improvable. Indeed, it is
known still in the one-dimensional case $ \ d = 1; \ $ for the multidimensional one it is sufficient to take as a trivial
{\it factorizable } function; say, when $ \ d = 2, \ $ one can choose

$$
f_0(x,y) := g_0(x) \ h_0(y), \ x \in X, \ y \in Y.
$$

\vspace{4mm}

 \section{Exponential bounds for distribution of positive multiple sums.}

\vspace{4mm}

 \ We intend to derive in this section the  uniform relative the  amount of summand $ \ |L| \ $ {\it exponential } bounds
for tail of  distribution of the r.v. $ \  S_L, \ $ based in turn on the moments bound obtained above as well as on the theory
of the so-called Grand Lebesgue Spaces (GLS). We recall now for readers convenience some facts about these spaces and supplement more.  \par

\vspace{4mm}

 \ These spaces are Banach functional space, are complete, and rearrangement
invariant in the classical sense, see [4], chapters 1, 2; and were investigated in particular in many works, see e.g. [5], [6]-[7], [13],
[19]-[20], [21], [22]-[25], and so one.\par

  \ They are closely related with the so-called {\it exponential} Orlicz spaces, see [5], [6], [7], [22], [23]-[25] etc. \par

\vspace{4mm}

  \ Denote for simplicity

$$
\nu_L(p) :=  W_L^{(d)}[f](p), \ \nu(p) :=  \sup_{L: |L| \ge 1} \psi_L(p), \eqno(6.1)
$$
and suppose

$$
 \exists b = \const \in (1, \infty]; \ \forall p \in (1,b) \ \Rightarrow \psi(p) < \infty. \eqno(6.2)
$$

\ Recall that the norm of the random variable $ \ \xi \ $ in the so-called Grand Lebesgue Space $ \ G \psi \ $ is defined as follows

$$
||\xi||G\psi \stackrel{def}{=} \sup_{p   \in (1,b)} \left\{ \frac{|\xi|_p}{\psi(p)}   \right\}. \eqno(6.3)
$$
 \ Here the {\it generating function} for these spaces $ \ \psi = \psi(p) \ $ will be presumed to be continuous  inside the open interval $ \ p \in (1,b) \ $
and such that

$$
\inf_{p \in (1,b) } \psi(p) > 0.
$$

  \ The inequalities (5.11) and (5.11a) may be rewritten as follows

$$
||S_L[f]||G\nu_L \le 1; \ \hspace{5mm} \sup_L  ||S_L[f]||G\nu \le 1. \eqno(6.4)
$$

 \ The so-called tail function $ \ T_{f}(y), \ y \ge 0 \ $ for arbitrary (measurable) numerical valued function (random variable, r.v.)
 $ \  f \ $ is defined as usually

$$
T_{f}(y) \stackrel{def}{=}  \max ( {\bf P}(f \ge y), \  {\bf P}(f \le -y) ), \ y \ge 0.
$$

 \ Obviously, if the r.v. $ \ f \ $ is non-negative, then

$$
T_{f}(y) =  {\bf P}(f \ge y), \  y \ge 0.
$$

 \ It is known that and if  $  \ f \in G\psi, \ ||f||G\psi = 1, \  $ then

$$
T_{f}(y) \le  \exp \left( -\zeta_{\psi}^*(\ln(y) \right),  \  y \ge e \ \eqno(6.5)
$$
where

$$
\zeta(p) = \zeta_{\psi}(p) := p \ \ln \psi(p).
$$

 \ Here the operator (non-linear) $ \ f \to f^* \ $   will denote the famous Young-Fenchel, or Legendre transform

$$
f^*(u) \stackrel{def}{=} \sup_{x \in \Dom(f)} (x \ u - f(x)).
$$

\vspace{5mm}

 \ We deduce by means of theorem 5.1 and property (6.5) \par

\vspace{5mm}

 \ {\bf Proposition 6.1.}

$$
T_{S_L[f]}(y) \le  \exp \left(  \ - \nu_L \left( \ \ln(y) \ \right) \ \right),  \  y \ge e;  \eqno(6.6)
$$

$$
\sup_L \ T_{S_L[f]}(y) \le  \exp \left( \ -\nu \left( \ \ln(y) \ \right) \ \right),  \  y \ge e. \eqno(6.6a)
$$

\vspace{5mm}

 \ {\bf Example 6.1}. \par
 \ Let us   bring  an example, see [30] for the centered r.v.  Let $  \  m = \const > 1 \  $ and define $  \  q = m' = m/( m-1). \  $
Let also $  \  R = R(y), \ y > 0 \ $ be positive  continuous differentiable {\it  slowly varying  } at infinity function such that

$$
\lim_{\lambda \to \infty} \frac{R(y/R(y))}{R(y)} = 1. \eqno(6.7)
$$
 \ Introduce a following $ \ \psi \ - \ $ function

$$
\psi_{m,R} (p) \stackrel{def}{=} p^{1/m} R^{-1/(m-1)} \left(  p^{ (m-1)^2/m  }  \right\}, \ p \ge 1, m = \const > 1, \eqno(6.7a)
$$

 \ Suppose

$$
\nu(p) \le \psi_{m,R} (p), \ p \in [1, \infty);
$$
then [19]-[20], [30] the correspondent exponential tail function has a form

$$
T^{(m,R)}(y) \stackrel{def}{=} \exp \left\{  - C(m,R) \ \ y^m \ R^{ m-1} \left(y^{m-1} \right)  \right\}, \ C(m,R) > 0, \ y \ge 1; \eqno(6.7b)
$$
so that

$$
\sup_L T_{S_L}(y) \le T^{(m,R)}(y), \ y \ge 1. \eqno(6.8)
$$

\vspace{4mm}

 \ A particular cases: $ \  R(y) = \ln^r (y+e), \ r = \const, \ y \ge 0; $ then the correspondent generating  functions has a form (up to multiplicative constant)

$$
\psi_{m,r}(p) = \ p^{1/m}  \ \ln^{-r}(p), \ p \in [2, \infty),  \eqno(6.9a)
$$
and the correspondent tail function has a form
$$
T^{m,r}(y) = \exp \left\{ \ - K(m,r) \ y^m \ (\ln y)^{r}   \ \right\}, \ K(m,r) > 0, \ y \ge e. \eqno(6.9b)
$$

 \ Many other examples may be found  in [19], [20], [22], [30] etc. \par

 \vspace{4mm}

\ {\bf  Example 6.2. }   Let the function $  \ f:  X = \otimes_{s=1}^d X_s \to R \ $ be from the degenerate representation

$$
f(\vec{x}) = \sum_{\vec{k} \in Z_+^d, \ N(\vec{k}) \le M} \lambda(\vec{k}) \ \prod_{j=1}^d g^{(j)}_{k_j}(x_j), \eqno(5.2a)
$$
for some constant integer value $ \ M, \ $  finite or not, where all the functions $ \   g_k^{(j)}(\cdot) \  $ are in turn non - negative:
$ \   g_k^{(j)}(\xi(k)) \ge 0. \ $ Recall the  denotation: $ \ M = \deg[f]. $ \par

 \ {\it Suppose here  and in what follows in this section that}

$$
\sum_{\vec{k} \in Z_+^d, \ N(\vec{k}) \le M} | \ \lambda(\vec{k}) \ | \le 1 \eqno(6.10)
$$
and that each the non-negative r.v. $ \  g_k^{(j)}(\xi(k)) \  $ belongs to  some $ \  G\psi_k \ - \ $ space  uniformly relative the index $ \ j: \ $

$$
 \sup_j| \ g_k^{(j)}(\xi(k)) \ |_p \le \psi_k(p). \eqno(6.11)
$$

 \ Of course, as a capacity of these functions may be picked the natural functions for the r.v. $ \ g_k(\xi(k)): \  $

$$
\psi_k(p) \stackrel{def}{=}  \sup_j |g_k^{(j)}(\xi(k))|_p,
$$
if the last function is finite for some non-trivial interval $ \ [2, a(k)), \ $ where $ \ a(k) \in (2, \infty]. \ $ \par

 \ Obviously,

$$
|f(\vec{\xi})|_p \le \prod_{k=1}^d \psi_k(p),
$$
and the last inequality  is exact if for instance $ \ M = 1 \ $ and all the functions $ \ \psi_k(p) \ $ are natural  for the  family of the
r.v. $ \ g^{(j)}_k(\xi(k)). \  $\par

\vspace{4mm}

 \ Define the following $ \ \Psi \ -  \ $ function

$$
\beta(p) = \kappa_d[\vec{\xi}](p) \stackrel{def}{=}  \prod_{k=1}^d \psi_k(p).
$$

  \ The assertion of proposition (5.1) gives us the estimations

$$
\sup_{L: |L| \ge 1} ||S_L[f]||G\kappa \le 1  \eqno(6.12)
$$
and hence

$$
\sup_{L: |L| \ge 1} T_{S_L[f]}(u) \le \exp \left( -v^*_{\kappa}(\ln u) \right), \ u \ge e, \eqno(6.12b)
$$
with correspondent Orlicz norm estimate. \par

\vspace{4mm}

\ {\bf  Example 6.3. } \par

 \ Suppose again that

$$
\sum_{\vec{k} \in Z_+^d, \ N(\vec{k}) \le M} \ | \lambda(\vec{k}) \ | \le 1
$$
and that the  arbitrary r.v. $ \  g^{(j)}_k(\xi(k)) \  $ belongs uniformly relative the index $ \ j \ $ to the correspondent
$   \  G\psi_{m(k), \gamma(k)} \  $ space:

$$
\sup_j | \ g^{(j)}_k(\xi(k)) \ |_p \le p^{1/m(k)} \ [\ln \ p]^{\gamma(k)}, \ p \ge 2, \ m(k) > 0, \ \gamma(k) \in R,  \eqno(6.13)
$$
or equally

$$
\sup_j T_{g^{(j)}_k(\xi(k))}(u) \le \exp  \left( - C(k) \ u^{m(k)} \ [\ln u]^{  - \gamma(k) }    \right), \ u \ge e. \eqno(6.13a)
$$

 \ Define the following variables:

$$
m_0 := \left[ \sum_{k=1}^d 1/m(k)  \right]^{-1}, \ \gamma_0 := \sum_{k=1}^d \gamma(k) , \eqno(6.14)
$$

$$
\hat{S}_L = \hat{S}_L[f] := e^d \ C^{-d}_R \ S_L. \eqno(6.15)
$$

 \ We conclude  by means of the proposition 5.1

$$
\sup_{ L: |L| \ge 1} \left| \left|  \hat{S}_L \right| \right| G\psi_{m_0, \gamma_0} \le 1 \eqno(6.16)
$$
and therefore

$$
\sup_{ L: |L| \ge 1} T_{\hat{S}_L}(u) \le \exp \left\{  - C(d,m_0, \gamma_0) \ u^{m_0} \ (\ln u)^{ - \gamma_0}   \right\},
\ u > e. \eqno(6.17)
$$

\vspace{4mm}

{\bf Example  6.4.}

\vspace{4mm}

 \  Let us consider  as above the following $  \   \psi_{\beta}(p) \  $ function

$$
  \psi_{\beta,C}(p)  :=  \exp \left( C p^{\beta} \right), \  C, \ \beta = \const > 0, \ p \in [1,\infty). \eqno(6.18)
$$
see example 6.3, (6.15)-(6.17). \par

\ Let $ \  g_k^{(j)}(\xi(k)) \ $ be  non-negative independent random variables belonging to the certain $ \  G \psi_{\beta,C}(\cdot) \  $ space
uniformly relative the indexes $ \ k,j: $

$$
\sup_j \sup_k ||g_k^{(j)}(\xi(k))||  G \psi_{\beta,C} = 1, \eqno(6.19)
$$
or equally

$$
\sup_j \ \sup_k T_{g_k^{(j)}(\xi(k)}(y) \le \exp \left(  \ - C_1(C, \beta) \ [  \ln(1 + y)   ]^{1 +1/\beta}  \  \right),  \ y > 0.
$$

  \ Then

$$
\sup_{L: \ |L| \ge 1} T_{S_L}(y) \le \exp \left(  \ - C_2(C, \beta) \ [  \ln(1 + y)   ]^{1 +1/\beta}  \  \right),  \ y > 0, \eqno(6.20)
$$
or equally

$$
\sup_{L: |L| \ge 1} || S_L[f]||  G \psi_{\beta,C_3(C, \beta)} =  C_4(C, \beta) < \infty. \eqno(6.20a)
$$

\vspace{4mm}

\ {\bf Example 6.5.} Suppose now that the each non-negative random variable $ \  g_k^{(j)}(\xi(k))  \  $  belongs uniformly relative the index $ \ j \ $
to certain $ \  G\psi^{<b(k), \theta(k)>} \  $ space, where $ \  b(k) \in (2, \infty), \ \theta(k) \in R: \ $

$$
\sup_j || \ g_k^{(j)}(\xi(k)) \ || G\psi^{<b(k), \theta(k)>} < \infty,
$$
where by definition

$$
\psi^{<b(k),\theta(k)>}(p) \stackrel{def}{=} C_1(b(k),\theta) \ (b(k)-p)^{ -(\theta(k) + 1)/b(k) }, \ 1 \le p < b(k).
$$

 \ This case is more complicates than considered before. \par

 \ Note that if the r.v. $ \  \eta \  $ satisfies the inequality

$$
T_{\eta}(y) \le C \ y^{-b(k)} \ [\ln y]^{ \theta(k)}, \ y \ge e,
$$
then $ \eta \in G\psi^{<b(k),\theta(k)>}, $ see the example 6.2.\par

 \ One can assume without loss of generality

$$
b(1) \le b(2) \le b(3) \le \ldots b(d).
$$

 \ Denote

$$
\nu_k(p) := \psi^{ < b(k), \theta(k) >}(p), \ b(0):= \min_k b(k),
$$
so that $ \ b(0) = b(1) = $

$$
 b(2) = \ldots = b(k(0)) < b(k(0) + 1) \le \ldots \le b(d), \ 1 \le k(0) \le d;
$$

$$
\Theta :=  \sum_{k=1}^{k(0)} (\theta(k) + 1)/b(0),
$$

$$
\upsilon(p) = \upsilon_{\vec{\xi}}[f](p) \stackrel{def}{=} \prod_{l=1}^{k(0)} \nu_l(p) =
C \cdot \left[ \ b(0) - p \  \right]^{  - \Theta }, \ 2 \le p < b(0).
$$

 \ Obviously,

$$
\ \prod_{k=1}^d \ \nu_k(p) \le C(d) \ \upsilon(p) =  C  \ \left[ \ b(0) - p \  \right]^{  - \Theta },  \
 C = C_d(\vec{\xi}, \vec{b}, \vec{\theta}, k(0)).
$$

 \ Thus, we obtained under formulated above conditions

$$
\sup_{L: |L| \ge 1} |S_L|_p \le C_2 \ (b(0) - p)^{-\Theta}, \ p \in [2, b(0))
$$
with the correspondent tail estimate

$$
\sup_{L: |L| \ge 1} T_{S_L}(y) \le C_3 \ y^{-b(0)} \ [ \  \ln y  \ ]^{ \ b(0) \ \Theta}, \ y \ge e.
$$

\vspace{4mm}

 \section{Upper bounds for these statistics. }

  \vspace{4mm}

 \  {\bf A.}  A simple lower estimate in the Klesov's (3.4) inequality may has a form

$$
\sup_{L: |L| \ge 1} \left| S^{(2)}_L   \right|_p  \ge \left| S^{(2)}_1  \right|_p =
 \ |g(\xi)|_p \ |h(\eta)|_p, \ p \ge 2, \eqno(7.1)
$$
as long as  the r.v. $ \   g(\xi), \ h(\eta) \ $ are independent. \par

 \ Suppose now that $ \ g(\xi) \in G\psi_1 \ $ and $ \  h(\eta) \in G\psi_2, \  $ where $ \ \psi_j \in \Psi(b), \ b = \const \in (2, \infty]; \ $
for instance $ \ \psi_j, \ j = 1,2 \ $ must be the natural functions  for these  r.v. Put $ \  \nu(p) = \psi_1(p) \ \psi_2(p);  \  $  then

$$
   \nu(p) \le \sup_{L: |L| \ge 1} \left| S^{(2)}_L   \right|_p \le K_1^d \cdot  \ \nu(p), \ K_1< \infty. \eqno(7.2)
$$

 \ Assume in addition that $ \ b < \infty; \ $ then $ \ K_1 \le C(b) < \infty. \ $ We get to the following assertion. \par

\vspace{4mm}

{\bf Proposition 7.1.} We deduce under formulated above in this section conditions

$$
 1 \le  \frac{ \sup_{L: |L| \ge 1} \left|S_L \right|_p}{\nu(p)} \le C^d(b) < \infty, \ p \in [2,b). \eqno(7.3)
$$

\vspace{4mm}

  \ {\bf B. \  Tail approach.}  We will use the example  6.2 (and notations therein. ) Suppose in addition that all the (independent) r.v. $  \ \xi(k) \ $
have the following tail of distribution

$$
T_{g_l(\xi(k))}(y) = \exp \left(  \ -  [\ln(1 + y)]^{1 +1/\beta} \ \right), \ y \ge 0, \ \beta = \const > 0,
$$
i.e. an unbounded support. As we knew,

$$
\sup_{L: \ |L| \ge 1} T_{S_L}(y) \le \exp \left(  \ - C_5(\beta,d) \ [  \ln(1 + y)   ]^{1 +1/\beta}  \  \right),  \ y > 0,
$$
  On the other hand,

$$
\sup_{L: \ |L| \ge 1} T_{S_L}(y) \ge T_{S_1}(y) \ge  \exp \left(  \ - C_6(\beta,d) \  [\ln(1 + y)]^{1 +1/\beta} \ \right), \ y > 0. \eqno(7.4)
$$

\vspace{4mm}

 \ {\bf C. An example.} Suppose as in the example 6.1 that the independent centered r.v. $ \  g^{(j)}_k(\xi(k)) \ $ have the  standard Poisson
 distribution:$ \ \Law(\xi(k) ) = \Poisson(1), \ k = 1,2,\ldots,d. \ $ Assume also that in the representation (5.2a)  $ \ M = 1 \ $
(a limiting degenerate case). As long as

$$
|g_k^{(j)}(\xi(k))|_p  \asymp  C \frac{p}{\ln p}, \ p \ge 2,
$$
we conclude by virtue of theorem 5.1

$$
\sup_{L: |L| \ge 1} \left|  \ S_L \right|_p \le C_2^d \ \frac{p^{2d}}{ [\ln p]^{2d}}, \ p \ge 2, \eqno(7.5)
$$
therefore

$$
\sup_{L: |L| \ge 1} T_{S_L}(y) \le \exp \left(  - C_1(d) \ y^{1/(2d)} \ [\ln y]^{2d}  \right), \ y \ge e. \eqno(7.6)
$$

 \ On the other hand,

$$
\sup_{L: |L| \ge 1} \left| S_L  \right|_p  \ge \left|  S_1 \right|_p \ge C_3(d) \ \frac{p^d}{ [\ln p]^d},
$$
and following

$$
\sup_{L: |L| \ge 1} T_{S_L}(y) \ge \exp \left(  - C_4(d) \ y^{1/d} \ [\ln y]^{d}  \right), \ y \ge e. \eqno(7.7)
$$

\vspace{4mm}

 \section{Concluding remarks. }

\vspace{4mm}

\ {\bf A.} \ It is interest by our opinion to generalize obtained in this report results onto the
mixing sequences or onto martingales, as well as onto the multiple integrals instead sums. \par

\vspace{4mm}

{\bf B.} \ Perhaps, a more general results may be obtained by means of the so-called method of
majorizing measures, see [1]-[3], [11], [29], [31]-[35].  \par

\vspace{4mm}

{\bf C. } \ Possible applications: statistics and Monte-Carlo method, alike [8], [10] etc.\par

\vspace{4mm}

{\bf D.} \ It is interest perhaps  to generalize the assertions of our theorems onto the sequences of domains $ \ \{ \ L \ \} \ $
tending to ``infinity'' in the van Hove sense, in the spirit of an articles [26]-[27], [30]. \par

\vspace{4mm}

 {\bf References.}

 \vspace{4mm}

{\bf 1. Bednorz W.} (2006). {\it A theorem on Majorizing Measures.} Ann. Probab., {\bf  34,}  1771-1781. MR1825156.\par

 \vspace{4mm}

{\bf 2. Bednorz  W.} {\it  The majorizing measure approach to the sample boundedness.} \\
arXiv:1211.3898v1 [math.PR] 16 Nov 2012 \\

\vspace{4mm}

{\bf 3. Bednorz W.} (2010), {\it Majorizing measures on metric spaces.}
C.R. math. Acad. Sci. Paris, (2010), 348, no. 1-2, 75-78, MR2586748. \\

\vspace{4mm}

{\bf 4. Bennet C., Sharpley R. } {\it Interpolation of operators.} Orlando, Academic
Press Inc., (1988). \\

\vspace{4mm}

{\bf 5. Buldygin V.V., Kozachenko Yu.V.} {\it Metric Characterization of Random
Variables and Random Processes. } 1998, Translations of Mathematics Monograph,
AMS, v.188. \\

\vspace{4mm}

{\bf 6. A. Fiorenza.} {\it Duality and reflexivity in grand Lebesgue spaces.} Collect.
Math. 51, (2000), 131-148.\\

\vspace{4mm}

{\bf 7.  A. Fiorenza and G.E. Karadzhov.} {\it Grand and small Lebesgue spaces and
their analogs.} Consiglio Nationale Delle Ricerche, Instituto per le Applicazioni
del Calcoto Mauro Picone”, Sezione di Napoli, Rapporto tecnico 272/03, (2005). \\

\vspace{4mm}

 {\bf 8.  Frolov A.S., Chentzov  N.N.} {\it On the calculation by the Monte-Carlo method definite integrals depending on the parameters.}
 Journal of Computational Mathematics and Mathematical Physics, (1962), V. 2, Issue 4, p. 714-718 (in Russian).\\

\vspace{4mm}

{\bf 9.  Gine,E.,  R.Latala,  and  J.Zinn.} (2000). {\it Exponential and moment inequalities for U \ - \ statistics.}
Ann. Probab. 18 No. 4 (1990), 1656-1668. \\

\vspace{4mm}

{\bf 10. Grigorjeva M.L., Ostrovsky E.I.} {\it Calculation of Integrals on discontinuous Functions by means of depending trials method.}
Journal of Computational Mathematics and Mathematical Physics, (1996), V. 36, Issue 12, p. 28-39 (in Russian).\\

\vspace{4mm}

{\bf 11. Heinkel  B.} {\it Measures majorantes et le theoreme de la limite centrale dan le  space C(S).}
Z. Wahrscheinlichkeitstheory. verw. Geb., (1977). {\bf 38}, 339-351.

\vspace{4mm}

{\bf 12. R.Ibragimov and Sh.Sharakhmetov.} {\it  The Exact Constant in the Rosenthal Inequality for Random Variables with Mean Zero.}
 Theory Probab. Appl., 46(1), 127–132. \\

\vspace{4mm}

{\bf 13. T. Iwaniec and C. Sbordone.} {\it On the integrability of the Jacobian under
minimal hypotheses.} Arch. Rat.Mech. Anal., 119, (1992), 129-143.\\

\vspace{4mm}

{\bf 14. Oleg Klesov.} {\it A  limit theorem for sums of random variables indexed by multidimensional indices. }
Prob. Theory and Related Fields,  1981, {\bf 58,} \ (3),  389-396. \\

\vspace{4mm}

 {\bf 15. Oleg Klesov.} {\it A  limit theorem for multiple sums of identical distributed independent
 random variables. }  Journal of Soviet Mathematics, September 1987, V. 38 Issue 6 pp. 2321-2326.
 Prob. Theory and Related Fields,  1981, {\bf 58, \ (3),}  389-396. \\

\vspace{4mm}

{\bf 16. Oleg Klesov.} {\it   Limit Theorems for Multi-Indexes Sums of Random Variables. } Springer, 2014.\\

\vspace{4mm}

{\bf 17. O. Klesov.} {\it Limit theorems for multi-indexed sums of random variables.}
 Volume 71 of Probability  Theory  and  Stochastic  Modeling.  Springer Verlag, Heidelberg, 2014. \\

\vspace{4mm}

{\bf 18. Korolyuk V.S., Borovskikh Yu.V.} (1994). {\it Theory of U-Statistics.} Kluwner Verlag, Dordrecht,
(translated from Russian).\\

\vspace{4mm}

{\bf 19. Kozachenko Yu. V., Ostrovsky E.I.}  (1985). {\it The Banach Spaces of
random Variables of subgaussian Type.} Theory of Probab. and Math. Stat. (in
Russian). Kiev, KSU, 32, 43-57. \\

\vspace{4mm}

{\bf 20.  Kozachenko Yu.V., Ostrovsky E., Sirota L} {\it Relations between exponential tails, moments and
moment generating functions for random variables and vectors.} \\
arXiv:1701.01901v1  [math.FA]  8 Jan 2017 \\

\vspace{4mm}

{\bf 21. E. Liflyand, E. Ostrovsky and L. Sirota.} {\it Structural properties of Bilateral Grand Lebesgue Spaces. }
Turk. Journal of Math., 34, (2010), 207-219. TUBITAK, doi:10.3906/mat-0812-8 \\

\vspace{4mm}

{\bf 22. Ostrovsky E.I.} (1999). {\it Exponential estimations for Random Fields and
its applications,} (in Russian). Moscow-Obninsk, OINPE.\\

\vspace{4mm}

{\bf 23. Ostrovsky E. and Sirota L.}  {\it Sharp moment estimates for polynomial martingales.}\\
arXiv:1410.0739v1 [math.PR] 3 Oct 2014 \\

\vspace{4mm}

{\bf 24. Ostrovsky E. and Sirota L.}  {\it  Moment Banach spaces: theory and applications. }
HAIT Journal of Science and Engineering C, Volume 4, Issues 1-2, pp. 233-262.\\

\vspace{4mm}

{\bf 25. Ostrovsky E. and Sirota L.} {\it  Schl\''omilch and Bell series for Bessel's functions, with probabilistic applications.} \\
arXiv:0804.0089v1 [math.CV] 1 Apr 2008\\

\vspace{4mm}

{\bf 26. E. Ostrovsky, L.Sirota.} {\it  Sharp moment and exponential tail estimates for U-statistics.  }
arXiv:1602.00175v1  [math.ST] 31 Jan 2016 \\

\vspace{4mm}

{\bf 27. Ostrovsky E. and Sirota L.} {\it Uniform Limit Theorem and Tail Estimates for parametric U-Statistics.}
 arXiv:1608.03310v1 [math.ST] 10 Aug 2016\\

\vspace{4mm}

{\bf 28. Ostrovsky E.I.} {\it  Non-Central Banach space valued limit theorem and applications.} In: Problems of the theory of probabilistic
distributions.  1989, Nauka, Proceedings of the scientific seminars on Steklov's Institute, Leningrad, V.11, p. 114-119, (in Russian).\\

\vspace{4mm}

{\bf  29. Ostrovsky E. and Sirota L.} {\it Simplification of the majorizing measures method, with development.} \\
arXiv:1302.3202v1  [math.PR]  13 Feb 2013 \\

\vspace{4mm}

{\bf 30. Ostrovsky E., Sirota L.} {\it Non-asymptotic estimation for Bell function, with probabilistic applications.} \\
arXiv:1712.08804v1  [math.PR]  23 Dec 2017 \\

\vspace{4mm}

{\bf 31. Talagrand  M.} (1996). {\it Majorizing measure:  The generic chaining.}
Ann. Probab., {\bf 24,} 1049-1103. MR1825156.\\

\vspace{4mm}

{\bf 32. Talagrand M.} (2005). {\it The Generic Chaining. Upper and Lower Bounds of Stochastic Processes.}
Springer, Berlin. MR2133757.\\

\vspace{4mm}

{\bf 33. Talagrand M.}  (1987). {\it Regularity of Gaussian processes.}  Acta Math. 159 no. 1-2, {\bf 99,} \ 149, MR 0906527.\\

\vspace{4mm}

{\bf 34. Talagrand  M.} (1990).  {\it Sample boundedness of stochastic processes under increment conditions.}  \ Annals of Probability,
{\bf 18,} \ N. 1, 1-49, \ MR10439.\\

\vspace{4mm}

{\bf 35. Talagrand  M.} (1992). {\it A simple proof of the majorizing measure theorem.}
Geom. Funct. Anal. 2, no. 1, 118-125. MR 1143666.\\

\vspace{4mm}

\end{document}